\documentclass[sn-mathphys,Numbered]{sn-jnl}

\usepackage{graphicx}%
\usepackage{multirow}%
\usepackage{multicol}
\usepackage{amsmath,amssymb,amsfonts}%
\usepackage{amsthm}%
\usepackage{mathrsfs}%
\usepackage[title]{appendix}%
\usepackage{xcolor}%
\usepackage{textcomp}%
\usepackage{manyfoot}%
\usepackage{booktabs}%
\usepackage{algorithm}%
\usepackage{algorithmicx}%
\usepackage{algpseudocode}%
\usepackage{listings}%
\usepackage{cleveref}

\usepackage{tkz-euclide}
\usetikzlibrary{quotes,angles,shapes.geometric,decorations.pathreplacing,shadings,shadows}
\usepackage{tikz,pgf,pgfplots}

\xdefinecolor{darkblue}{rgb}{0,0.3,0.58}
\xdefinecolor{darkgreen}{rgb}{0,0.6,0.3}
\newcommand{\R}{\mathbb{R}}

\theoremstyle{thmstyleone}%
\theoremstyle{thmstyletwo}%

\theoremstyle{thmstylethree}%

\begin{document}

\title[Mixture of experts surrogate model for the homogenization of open-porous materials]{Mixture of experts surrogate model for the homogenization of open-porous materials}


\author*[1,2]{\fnm{Axel} \sur{Klawonn}}\email{axel.klawonn@uni-koeln.de}

\author*[1,2]{\fnm{Martin} \sur{Lanser}}\email{martin.lanser@uni-koeln.de}

\author*[1,3]{\fnm{Lucas} \sur{Mager}}\email{lucas.mager@uni-koeln.de}

\author*[3,4]{\fnm{Ameya} \sur{Rege}}\email{ameya.rege@utwente.nl}

\author*[1,2]{\fnm{Janine} \sur{Weber-Hamacher}}\email{janine.weber@uni-koeln.de}

\affil[1]{\orgdiv{Department of Mathematics and Computer Science}, \orgname{University of Cologne}, \orgaddress{\street{Weyertal 86-90}, \city{Cologne}, \postcode{50931}, \country{Germany}}}

\affil[2]{\orgdiv{Center for Data and Simulation Science}, \orgname{University of Cologne},  \orgaddress{\street{Albertus-Magnus-Platz}, \city{Cologne}, \postcode{50923}, \country{Germany}}}

\affil[3]{\orgdiv{Institute for Frontier Materials on Earth and in Space}, \orgname{German Aerospace Center}, \orgaddress{\street{Linder Höhe}, \city{Cologne}, \postcode{51147}, \country{Germany}}}

\affil[4]{\orgdiv{Department of Mechanics of Solids, Surfaces \& Systems}, \orgname{University of Twente}, \orgaddress{P.O. Box 217, Enschede, 7500 AE, \country{Netherlands}}}


\abstract{For open-porous materials, incorporating their microstructural properties into mechanical simulations poses a significant challenge for accurately capturing elastic deformation. 
 To deal with this difficulty, multiscale methods are a common tool to couple characteristics of the microstructure of the considered material with the macroscopic material behavior. However, when desiring a high accuracy, these multiscale computations can be computationally very expensive due to the large number of microscopic problems which need to be solved in each compute step.  Here, surrogate models that learn the mechanical response of the underlying constitutive model can significantly reduce the computational cost of multiscale approaches.  
 In previous work by some of the authors, beam frame models have been used to model the microstructure of open-porous materials which have been combined with neural network-based surrogate models to approximate the material behavior of a given RVE (repesentative volume element). In this work, we extend our previous study by training a more complex neural network model 
 to predict the mechanical behavior of several RVEs, differing in their maximum pore size and pore-size distribution. 
 Concretely, we focus on mixture of expert (MoE) models and compare different MoE architectures as well as their performance across different RVEs. This novel approach reduces the computational cost of simulating multiple RVEs as the MoE model does not require additional training when new RVEs are considered. }




\maketitle
\newpage

\section{Introduction}\label{secIntroduction}

 Owing to their unique combination of low density, high porosity, thermal and acoustic insulation, energy absorption and damping capacity, as well as their potential for functionalization, open-porous materials such as aerogels are of considerable interest for applications ranging from lightweight structures and extreme-temperature environments to biomedical technologies. 
Depending on the base material used for synthesis and the conditions for the sol-gel process or the drying, the typical open-porous nanostructure of aerogel materials which significantly dictates the material properties can vary a lot between different materials. In addition to the insulation properties, the mechanical behavior is also determined by the nanoscopic pore structure of the material. Therefore, for simulation methods which compute the mechanical response of aerogels, it is important to consider the respective nanostructure of the material. For numerical simulation methods, this poses a significant challenge. 

Regarding the micromechanical behavior, the research on the deformation of open-porous structures dates back to work of Gent and Thomas of 1959 \cite{gent-thomas-1959} which considers the deformation of foam structures. By modeling the fibrillar pore structure with Euler-Bernoulli beam elements, Gibson and Ashby introduced a more robust modeling method \cite{gibson-prsa-1982,gibson-ashby-prsa-1982}. The strain energy approach by Dement'ev and Tarakanov \cite{dementev-tarakanov-1970} was introduced to study the plastic deformation of foams. Based on this approach, Rege et al. provided a micromechanical constitutive model for the mechanical behavior of cellular solids \cite{rege-pre-2021}. Further refinement of this approach for the study of biopolymer aerogels yields good validation results regarding the mechanical behavior under large strains \cite{rege-mat-2021}.

For the analysis of macroscopic deformation behavior, multiscale methods are often the preferred approach for the simulation of materials with a heterogeneous micro- or nanostructure. With the application of multiscale models, it is possible to couple the macroscopic deformation response with the micromechanical behavior. Especially, the FE$^2$ method \cite{MieheBound,SMIT_Hom_1998,HacklSchroeder,Kouz_Hom_2001,FEYEL_FE2_1999} has established as a well-known approach for these applications. In general, the method is computationally very demanding due to many microscopic finite element problems which need to be solved on representative volume elements (RVEs). Additionally, these problems cannot be too small, since the RVE has to capture the effective properties of the material. To reduce the computational effort of solving the microscopic problems, machine learning based surrogate models are often developed to replace the computationally expensive FEM solver on the microlevel \cite{surr2,surr3,surr4,surr5,surr6,surr7,surr8,surr9,surr10,surr11,wang2024}.

The present paper builds upon the work of \cite{AerogelFE2} in which a multiscale method for open-porous materials is introduced. The approach considers a beam frame problem on the microscopic scale. By developing a neural network (NN) trained to predict the homogenized stress components resulting from the deformed RVE, the resulting surrogate model can be used to replace the beam frame model on the microscopic scale. In this method, the NN is trained to replicate the homogenized stress tensor resulting from the deformation of a single RVE. Our previous study~\cite{AerogelFE2} shows that the NN approach yields significantly reduced computation times and that the error which is introduced by the surrogate model stays relatively low.  However, as a limitation, the surrogate model presented in~\cite{AerogelFE2} lacks generalizability, as it is trained and tailored exclusively to one specific RVE and, consequently, a unique material composition. To overcome this limitation, we extend the given multiscale method from~\cite{AerogelFE2}, in the present work, and introduce a surrogate approach based on the mixture of experts (MoE) architecture. 
We introduce this model with the objective of developing a general-purpose surrogate model which takes a numerical representation of the considered nanostructure as an additional input. By considering the pore structure of the material, the MoE can be trained on the homogenized stress values of different open-porous materials.
 
  To the best of our knowlege, the MoE approach was first introduced in the field of machine learning (ML)
 by Jacobs, Jordan, Nowlan, and Hinton \cite{jordanjacobs1991MoE}; see also Jordan and Jacobs \cite{jordanjacobs1994MoE}. A notable MoE application is the multilayer MoE model introduced by Eigen \cite{eigen2013deepMoE} for an image classification problem based on the MNIST dataset which laid foundation for various further developments regarding ML with the MoE architecture. In the research of learned surrogate models for mechanical simulations, the MoE approach is often used for increasing the accuracy in simulation setups with complex input-output correlations \cite{StillerGatedPINNs2020, DanMoEPINNs2025, NabianMoESurrogateAerodynamics2025, BischofMoEPINNs2022}. For these applications, the MoE architecture can be advantageous compared to a single large NN and it can be directly combined with different approaches like physics-informed loss terms \cite{StillerGatedPINNs2020, DanMoEPINNs2025, BischofMoEPINNs2022}.

In the scope of the present work, we apply the MoE approach to the FE$^2$ method introduced in \cite{AerogelFE2} with the objective to develop a learned surrogate model capable of differentiating between different microscopic problems. For training the MoE surrogate model, we consider two training phases. In the first phase, each expert NN is trained on a separate dataset of deformation and stress values given from the computations for a single RVE. Consequently, each expert is specialized on a specific material and the expert NNs are comparable to the surrogate model introduced in~\cite{AerogelFE2}. Following the insights of \cite{AerogelFE2}, it is expected that each expert NN is capable of replicating the mechanical behavior of the material it is trained for. In the second training phase, a gating NN is trained to predict the best combination of the expert outputs based on given RVE inputs. For a given open-porous material and deformation gradient, the homogenized stress predicted by the MoE model is then computed from the weighted sum over the expert outputs.

The material presented in this article was previously included in Chapter 4, esp. Section 4.3.2, of the doctoral thesis of the third author \cite{MagerDiss2026}. The thesis explicitly identifies this material as part of joint work with the present coauthors and refers to the present article as work in preparation; see reference [83] in \cite{MagerDiss2026}.

\section{Computational homogenization approach for beam frame microstructure}\label{secFE2}
The multiscale approach which we apply for simulating the mechanical behavior of the open-porous material with respect to the given nanostructure has been introduced in~\cite{AerogelFE2}. The proposed approach is based on the FE$^2$ method~\cite{MieheBound,SMIT_Hom_1998,HacklSchroeder,Kouz_Hom_2001,FEYEL_FE2_1999}. It considers a finite element problem on the macroscopic level which does not resolve the pores on the micro- or nano-scale. At each integration point of the macroscopic finite elements a localized microscopic problem is attached. For these microscopic problems, we consider an RVE which is based on the nanostructure of the material and represents its most important features. 
 Given the three-dimensional connected fibrillar network that makes up the nanostructure of the material, the RVE is modeled as a beam frame structure. Homogenization over the RVE is done following the formulas introduced in~\cite{AerogelFE2}.

\subsection{Microscopic beam frame model}\label{secBF}
 Since we consider open-porous materials, constructing the geometry of the RVEs which reflects the most important features is crucial. Therefore, for creating an RVE for the microscopic problem, we apply the method presented in~\cite{DLR_biopoly_aerogel2019,DLR_biopoly_aerogel2021}. This approach focuses on the distribution of pore sizes and the porous fraction of the material. These properties are usually determined by experimental analyses from a given sample. The model framework yields a beam frame model which is expected to have similar structural properties as the nanostructure of the given material. The method ensures that the resulting RVE has periodic boundaries. Computing the solution of a three-dimensional microscopic beam frame problem with $n_{B}$ beam elements connected at $n_{N}$ vertices requires solving a linear system of equations
\begin{equation*}
    K\, u = F
\end{equation*}
for the deformation vector $u\in R^{6\cdot n_{N}}$ with the stiffness matrix $K\in R^{6\cdot n_{N}\times6\cdot n_{N}}$. The matrix $K$ is assembled from submatrices $K_e$ for each beam element $e$. The exact formulation of the submatrices based on~\cite{kassimali2011} is presented in~\cite{AerogelFE2}. The solution $u$ of the beam frame model features information about the deformation and the rotation of each node within the RVE. In three dimensions, this leads to 6 degrees of freedom (dof) for each node. 

\subsection{Homogenization approach}
\begin{figure}
\centering
\includegraphics[width=0.75\textwidth]{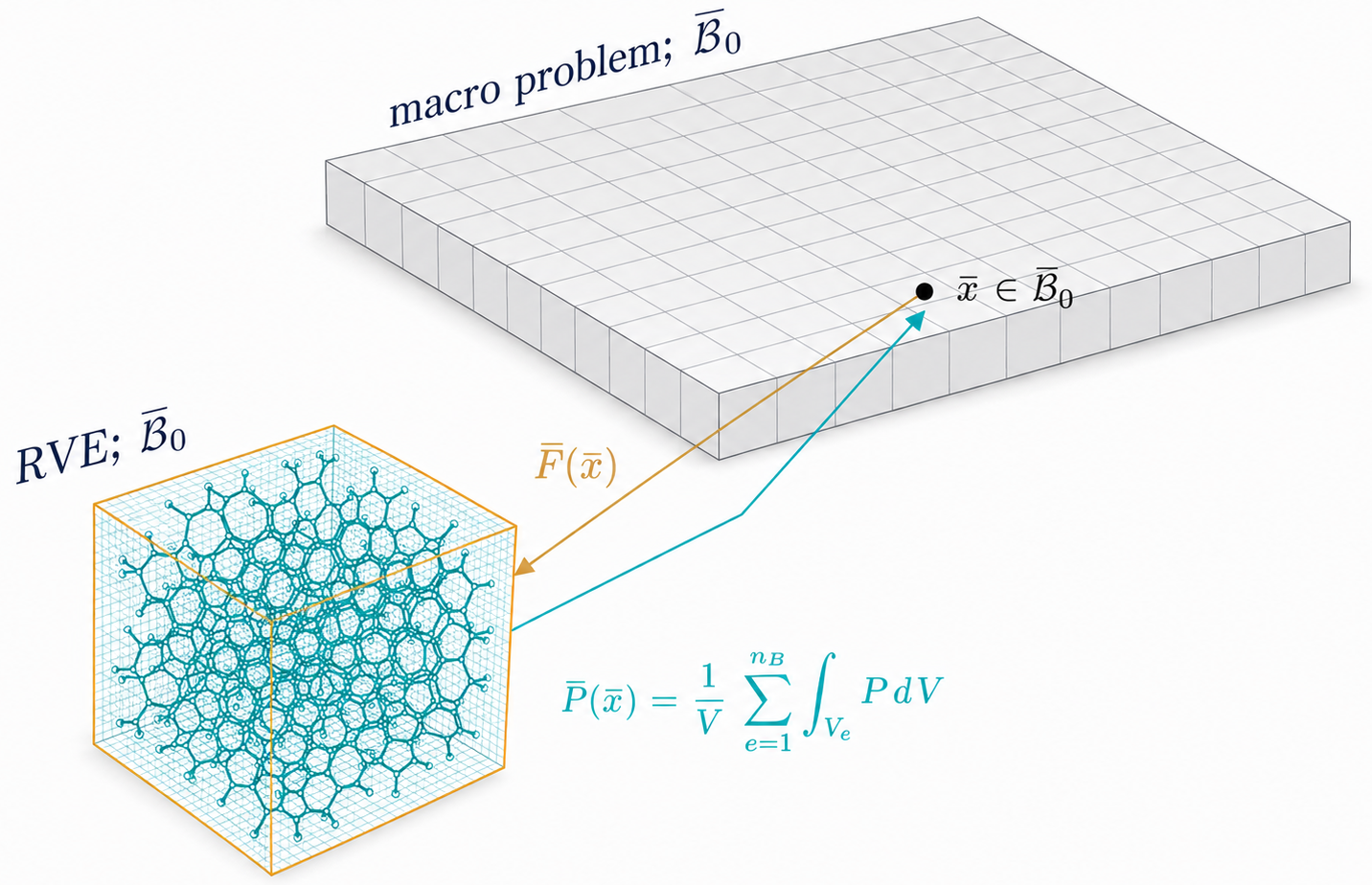}	
\caption{Visualization of the homogenization approach; the macroscopic problem is discretized with finite elements; in each Gauss integration point an RVE is attached (we show only one exemplary one here); the RVE is discretized using beam elements.}
\label{fig:fe2}
\end{figure}

For the macroscopic problem, we consider the momentum balance equation which can be presented in the weak form on the macroscopic reference domain $\overline{B}_0$ by the equation
\begin{align*}
	\int_{\overline{B}_0} \delta\bar{x}\left({\rm Div}_{\bar{x}}\overline{P}(\overline{F})-\Bar{f}\right)d\Bar{x}=0
\end{align*}
with the volume forces given by $\Bar{f}$. For the extrapolation of the integral, the macroscopic first Piola-Kirchhoff stress tensor $\overline{P}$ is computed for every Gauss quadrature point in $\overline{B}_0$ from the solution of the corresponding microscopic RVE problem. The boundary condition of a specific microscopic problem is based on the macroscopic deformation in  the corresponding quadrature point given by the macroscopic deformation gradient $\overline{F}$. More specifically,  Dirichlet boundary conditions are applied to the eight corners of the RVE and periodic boundary conditions are applied to the other boundary nodes of the beam frame model. With the resulting deformation and rotation vectors, it is possible to compute the average first Piola-Kirchhoff stress tensor for each beam element in the RVE. The computation of the average stress within one beam element is based on~\cite{ARL_stress_av}. The homogenized stress for the complete RVE is then computed as the sum over the Piola-Kirchhoff stress tensors of all beam elements
\begin{equation*}
    \overline{P}=\sum_{e=1}^{n_{B}} \int_{V_e} P dV
\end{equation*}
with $V_e$ being the volume of the beam element $e$. The relation between the macroscopic and one specific RVE problem is illustrated in \cref{fig:fe2}.

\subsection{Neural network-based surrogate model}
The described homogenization approach that considers the beam frame model on the microscopic scale is very versatile since it is easily applicable to many different aerogel materials with different nanostructures. Although the linear beam frame model can be evaluated faster than a similar full finite element discretization, it can still be quite computationally expensive to run a complete multiscale simulation with a complex microstructure. Particularly large macroscopic problems with many finite elements require high computing power. 

As described in~\cite{AerogelFE2} it is possible to replace the microscopic beam frame problem and the homogenization of the stresses with a neural network-based surrogate model which is trained to predict the mechanical behavior of a concrete open-porous material with a given nanostructure. We have previously shown that the surrogate approach is capable of yielding the same convergence behavior as the beam frame approach and the difference between the results which is caused by the  approximation error of the neural network (NN), stays small. These NNs take the macroscopic deformation gradient $\overline{F}$ as an input and directly predict the homogenized first Piola-Kirchhoff stress tensor $\overline{P}$ for the RVE. Since relatively small ML models are sufficient for predicting this relationship between deformation and stress, this method allows to run simulations with a significantly lower computational effort.  

 One disadvantage of the surrogate approach described in~\cite{AerogelFE2} is that the NN replacing the microscopic beam frame model is only trained to predict the stresses for one given material, that is, for one specific RVE. Using the surrogate approach for the multiscale simulation of a different material with a different nanostructure requires the training of a new NN.  Training a surrogate model creates significant overhead, requiring the computation of about $100\,000$ data samples in addition to the computations during training. The training of a NN for a single RVE takes approximately 30 minutes with running the computations on a Nvidia DGX-1 GPU. 
In the following section we introduce a more complex surrogate model based on the MoE approach which takes multiple nanostructures into account and is expected to yield a better versatility with regard to different RVEs and generalizes better to RVEs and thus nanostructures not included in the training data.

\section{Mixture of experts surrogate model}\label{secMoE}

 As explained above, the aim of this work is to develop a surrogate model on the microlevel that performs well for different types of material, that is, different RVEs. However, training a single neural network model for this purpose potentially requires a complex network architecture with a high number of parameters as well as a large amount of training data. In  preliminary tests, single NN architectures trained to predict different materials have also been strongly affected by overfitting. Therefore, the use of a single NN model has shown to be not beneficial for this application. 
Hence, the present approach makes use of the idea of conditional computing and training multiple expert models which are each specialized for a concrete task, that is, one specific material in form of a specific RVE; see, for example,~\cite{bengio2015conditional,shazeer2017outrageously}. Typically, these expert models are then aggregated by a linear or nonlinear gating network, depending on the networks' input and resulting in a mixture of experts (MoE) model. 
The early development of the mixture of experts approach has been based on the divide-and-conquer principle~\cite{Smith1985DnC}, which partitions a given problem space into multiple subspaces. Similar to this approach, in an  MoE model, each expert is trained to predict the solution within a given input subspace. The gating network is expected to identify the weight for each expert to compute the output as a weighted combination of all expert predictions. The combination is not fixed but automatically computed individually for each given input. 
  This decomposition is particularly attractive for the present application, as it allows the individual experts to learn the constitutive response of specific RVEs without requiring a single network to represent the entire heterogeneous response space. At the same time, the gating network enables the specialized expert predictions to be combined adaptively, providing a more flexible and potentially more generalizable surrogate across different RVEs. 

\begin{figure}
\centering
\includegraphics[width=1.0\textwidth]{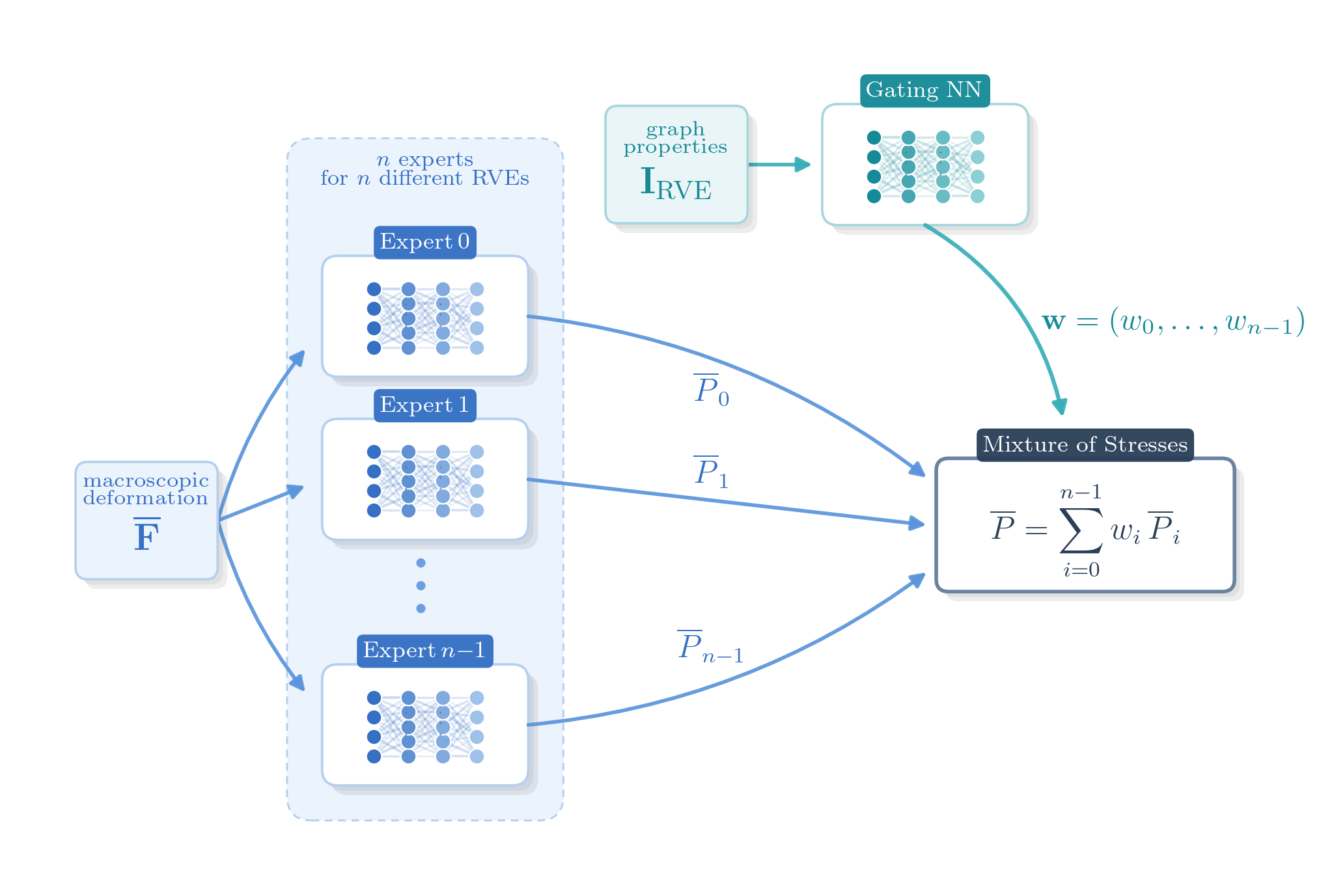}
    \caption{Illustration of the mixture of experts model. The input $\mathbf{\overline{F}}$ of the expert models refers to the macroscopic deformation gradient and the input $\mathbf{I_\text{\tiny RVE}}$ of the gating model refers to the parameterized representation of the input RVE.}
    \label{MoE:fig}
\end{figure}

 In the following, we define in detail our proposed MoE model that should replace the microscopic problem within the FE$^2$ framework. 
 In particular, our MoE model can be split into two parts which are also trained in different training phases. The structure of the MoE model is illustrated in \Cref{MoE:fig}. The first part of the full model which is presented in blue color  considers all individual expert models. These are NN models which are each pre-trained to predict the average stress components resulting from given deformations of a fixed RVE. So each expert model is related to one given nanostructure and is basically a NN-based surrogate model as introduced in our earlier work~\cite{AerogelFE2}. The expert models all take the same macroscopic deformation gradient ${\overline{F}}$ as an input. The output of the $i$-th expert is the resulting first Piola-Kirchhoff stress tensor ${\overline{P}_i}$.

The second part of the MoE model consists of the gating model which is expected to predict the importance of the expert models for a given nanostructure. It is not required that the nanostructure is included in the set of $n$ different nanostructures covered by the $n$ expert models. Hence, the input of the gating model is an aggregated representation of the nanostructure of the material which should be simulated. This representation of the nanostructure is notated as $\bf{I}_{RVE}$ in \Cref{MoE:fig}. We choose some parametrized representation of the structure and provide all details on $\bf{I}_{RVE}$ in~\Cref{secParam}. The output of the gating model is a vector with the dimension matching the number of expert models. This output vector $w$ defines the weights used for the weighted sum over all expert outputs to compute the final output of the full  MoE model. The combined output of the MoE model is then given by
\begin{equation*}
    \overline{P}=\sum_{i=0}^{n-1} w_i \cdot \overline{P}_i.
\end{equation*}
\subsection{Expert neural networks}\label{secExpertNN}

 In our approach, each of the expert models is assigned to a certain RVE and is trained to predict the resulting first Piola-Kirchhoff stress components resulting from the given deformation of the RVE. For the presented results in~\Cref{secTrain}, we use the same dense feedforward architecture for all expert models which is based on the results of \cite{AerogelFE2} and we also test different numbers of experts. For the MoE approach, it is generally not required to have the same architecture for each expert. 
 Here, we decided to use the same architecture for each expert, since all expert NNs are trained on very similar problems. We expect that any performance improvements resulting from the selection of varying expert architectures would be minimal at most.  
 In our case, each expert NN consists of three hidden layers with gelu activation and 256 neurons per layer. In addition to the fact that the suitability of the gelu activation function has been demonstrated in \cite{AerogelFE2}, its differentiability and smoothness make it particularly well-suited for approximating the continuous nonlinear stress-strain response of the underlying constitutive model.  All expert models take the deformation gradient as an input which consists of nine features in the three-dimensional case. For computing the nine output components representing the first Piola-Kirchhoff stress tensor, a linear activation is used in the final layer. 

\begin{figure}
    \centering
    \includegraphics[height=6cm, trim={0cm 0cm 0cm 0cm}, clip]{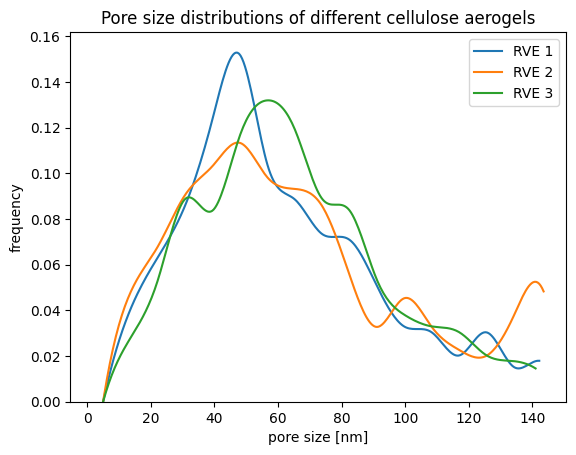}
    \caption{Density plots of the pore-size distributions for different randomly generated aeorgel RVEs. Image and data have first been published in \cite{MagerDiss2026}.}
    \label{Distr:fig}
\end{figure}

For generating the training data for one expert model we first need to create an RVE. The pore-size distribution and the porous fraction for the RVE are determined randomly but are set to be similar to values obtained from experimental analyses of cellulose aerogels \cite{DLR_biopoly_aerogel2018}. Three examples of pore-size distributions are shown in \Cref{Distr:fig}.
The data used for training the expert models is generated by computing the homogenized stresses of the beam frame problem for the specific RVE. We consider a dataset of 100,000 randomly chosen deformation gradients which is split into 90,000 training data samples and 10,000 validation data samples. The size of the training data is chosen to achieve a large variety of inputs and allows each model to generalize for different deformations. 

\begin{figure}
    \begin{minipage}{0.32\textwidth}
        \includegraphics[width=\textwidth]{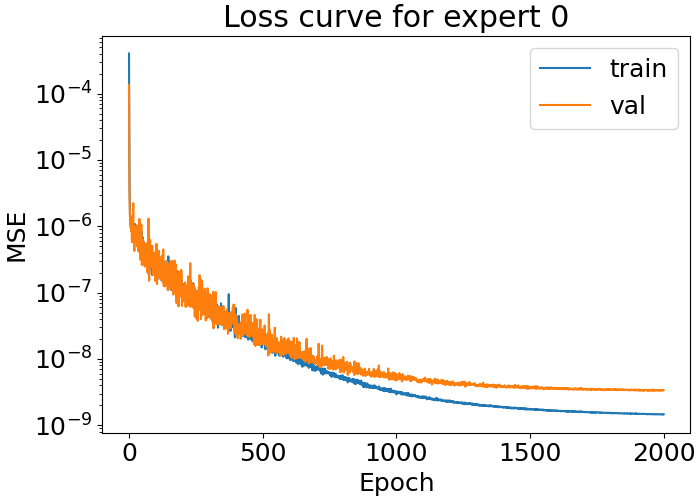}
    \end{minipage}
    \begin{minipage}{0.32\textwidth}
        \includegraphics[width=\textwidth]{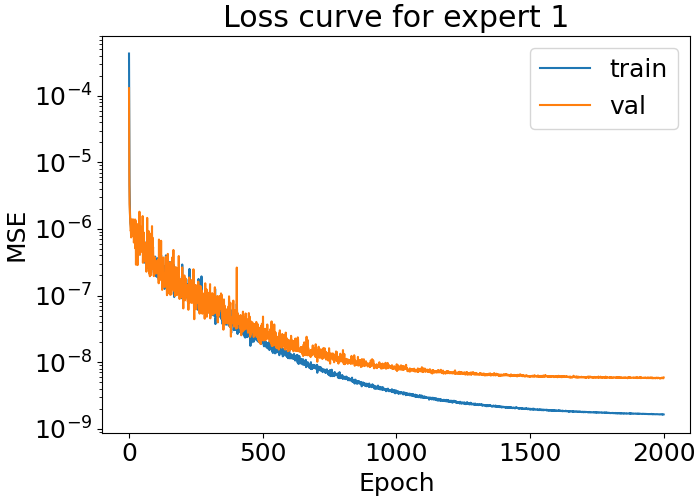}
    \end{minipage}
    \begin{minipage}{0.32\textwidth}
        \includegraphics[width=\textwidth]{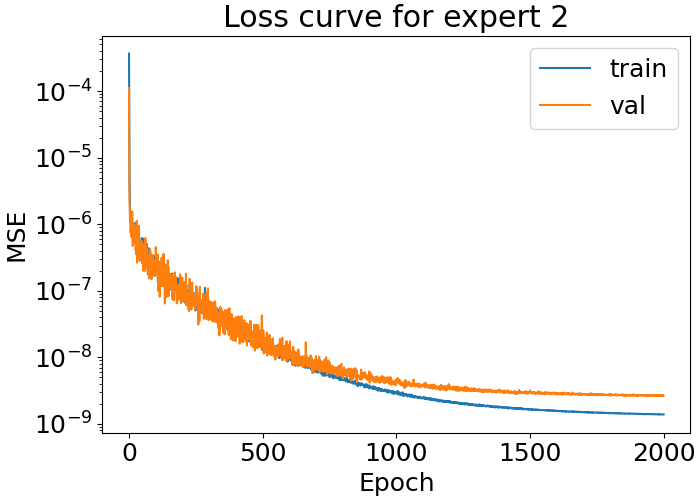}
    \end{minipage}
    \caption{Loss curves for the training and validation of the first three expert models. The training loss is shown in blue and the orange curve shows the behavior of the validation loss over $2\,000$ epochs. Images and data have first been published in \cite{MagerDiss2026}.}
     \label{expert_loss:fig}
\end{figure}

Before the weights of the gating model are adjusted, each of the expert networks is trained individually on its generated data  in a first phase. These training processes can, in principle, be computed in parallel since all expert models are fully independent from each other. 
 Each expert is trained for 2,000 epochs with an adaptive learning rate scheduler. The loss curves for the training of the first three experts are presented in \Cref{expert_loss:fig}. The diagrams show that the training and validation loss reduce significantly over the full number of epochs. The validation loss only ends up slightly higher than the training loss which indicates that the models should not be too much affected by overfitting. The same training behavior is also observed with the other expert models. 

\subsection{Parameterized representation of the RVE}\label{secParam}
As described in \Cref{secBF}, the method for creating the beam frame RVE takes the pore-size distribution and the porous fraction of the considered material into account. Therefore, we decided to use these properties also for the input $\bf{I}_\text{RVE}$ of the gating model. Additionally, we add some more informative statistics as input data, as, for example, orientations of the beams in the RVE. The complete set of input features is given by the following set of values:
\begin{itemize}
    \item[1)] {Mean} of pore-size distribution
    \item[2)] {Median} of pore-size distribution
    \item[3)] {Standard deviation} of pore-size distribution
    \item[4)] {Skewness} of pore-size distribution
    \item[5)] {Kurtosis} of pore-size distribution
    \item[6)] {Normalized orientation} of the beam elements in $x$-direction
    \item[7)] {Normalized orientation} of the beam elements in $y$-direction
    \item[8)] {Normalized orientation} of the beam elements in $z$-direction
    \item[9)] {Porous fraction}
    \item[10)] {Mean degree} of the nodes within the RVE
\end{itemize}
In total, the vector $\bf{I}_\text{RVE}$ contains ten values and for the remainder of this work we denote the number of input features for the gating model as $n_{in}=10$. The first five statistical values describe the pore-size distribution while the next three values $O_x,\, O_y,\, O_z$ describe the orientation of the beam elements within the coordinate system.  They are computed as the sum over the absolute orientation vectors weighted by their cross-sectional area, that is,
\begin{equation*}
    \left(\begin{array}{c}
        O_x \\
        O_y \\
        O_z
    \end{array}\right) = \sum_{e = 1}^{n_B}(r_e^2\cdot \pi)\cdot \lvert v_{\rm{end}(e)} - v_{\rm{start}(e)}\rvert,
\end{equation*}
with $v_{\rm{start}(e)}$ and $v_{\rm{end}(e)}$ being the coordinates of the start and end nodes of a beam element $e$. Here, the radius of the beam element is denoted by $r_e>0$. The orientation variables are normalized in a way that their mean equals  one. The porous fraction is given from the experimental analyses of the material which were previously mentioned in \Cref{secBF} and the mean degree variable is simply obtained by $2\frac{n_B}{n_N}$. 

\subsection{Gating network}
 In our MoE approach, the gating network is basically trained to solve a classification problem since it predicts the probability distribution over all experts to be picked for the simulation of a given nanostructure. The application of NNs for solving classification problems is very well studied and especially the \textit{softmax} activation function has been established as a good choice for classification networks and has already been used in such applications in the early 1990s~\cite{Bridle1990}.

For our gating network, we again consider a feedforward model architecture. The network consists of two hidden layers and one output layer with the number of output neurons matching the number of experts. For the output activation function, we consider the \textit{softmax} function since it provides probability values which sum up to 1. 

 As already mentioned, the parameters (weights and biases) of the gating network are learned in a second training phase. In general, in this second training phase, it is possible to either fix the weights and biases of all expert networks or to keep the expert models also trainable. In the latter case, the MoE model is trained in a coherent approach and all the experts, pre-trained in the first training phase, can be further optimized in the second training phase. This is not the case in the first approach, where the experts do not change during the training of the gating network. We compare the performance of the two training methods for the gating network in \Cref{secTrain}.

\section{Training and validation for the mixture of experts model}\label{secTrain}
 This section assesses the MoE model on the microscopic scale for various RVEs. It compares the MoE model's performance with the beam frame approach and models trained for single RVEs. We consider different gating model architectures with different numbers of experts and analyse the relevance of the different statistical input features representing each RVE.

\subsection{Analysis of loss curves and training progress}
 As described in \Cref{secExpertNN}, we use a two stage training process to optimize the weights and bias values of the MoE model with respect to the training and validation data.  Hence, in a first phase, we exclusively train the expert models, whereas subsequently, in a second phase, the gating model is trained. 
This can either be done with the expert weights and biases set as trainable which allows adjusting all weights of the MoE model  simultaneously, or by freezing the expert weights to ensure that the training process of the first training stage is retained,  
 that is, each expert is able to predict the stresses for a given RVE with high precision.

 {\bf First phase: training of expert NNs.} We create 250 different RVEs and computed, using the beam frame model, the average stress $\overline{P}$ for $100\,000$ different macroscopic deformations. In total, this results in a dataset of 25 million pairs of macroscopic deformations $\overline{F}$ and averaged stresses $\overline{P}$, or, in other words, a set of $100\,000$ pairs per RVE. With these 250 sets of data, we train 250 expert NNs using the optimal model architecture from~\cite{AerogelFE2}. We split the datasets into 90\% training and 10\% validation data for the training of each expert. The training of the 250 NNs is independent of each other and there is no overlap in the training datasets, since they are all specialized on exactly one RVE.

 {\bf Second phase: training of the MoE.}  
 As mentioned above, in the second training stage, either only the parameters of the gating network are optimized while the expert networks remain frozen, or all parameters of the MoE model are jointly fine-tuned.  
In both cases, we consider a training dataset of $1\,800\,000$ samples, each consisting of an input deformation $\overline{F}$ for the experts, an input $\bf{I}_\text{RVE}$ for the gating NN, and an output $\overline{P}$ of the global MoE. The 1.8 million data samples were randomly selected from the dataset of $25\,000\,000$ samples that was already used for the training of the experts in the first phase. Let us note that the statistical features representing the specific RVE $\bf{I}_\text{RVE}$  have not been used to train the experts and only serve as an input of the gating network; see~\Cref{secParam} for a description of $\bf{I}_\text{RVE}$.

 To control the generalization properties of our model during the training process, we consider two different validation datasets. This is necessary since the final MoE should be able to generalize across different deformations and RVEs.  The first validation dataset consists of $200\,000$ samples that are randomly selected from the same large dataset of $25\,000\,000$ samples as the training dataset. The method for the random selection ensures that the validation data does not share any samples with the training data. 
 Thus, this validation dataset probably differs from the training dataset with regard to the deformation inputs but not with regard to the inputs of the gating model. This dataset is used to evaluate the model's performance to generalize with respect to the deformation input. The second set of validation data consists of $400\,000$ samples which were randomly selected from a different, separate dataset containing $5\,000\,000$ data points. This dataset is built from $100\,000$ deformation gradients and resulting stress components computed with 50 RVEs different from the 250 RVEs used for the generation of the training data.   Thus, the second set of validation data differs from the set of training data not only with regard to the deformation inputs but also with regard to the parameterized RVE representation.  This dataset enables the evaluation of the models ability to generalize with respect to different RVE inputs and to evaluate the ability of the gating NN to learn from the input $\bf{I}_\text{RVE}$. To differentiate between the two validation datasets and to simplify the naming, we refer to the first validation dataset as $val_{\rm F}$ and to the second as $val_{\rm RVE}$.

Monitoring the loss of the MoE models during training regarding the training data as well as the validation datasets is crucial since the dimension of these loss values indicate the model's ability to replicate the behavior of the beam frame model on the microscopic scale of a material. As a consequence, for evaluating a model's performance, we first compare the loss curves on the validation datasets.

 Additionally, for comparing different MoE architectures, we also compute an entropy based confidence measurement for the respective gating model NN$^\text{gate}$. For the computation, we use the entropy formulation
\begin{equation*}
    H = -\sum_{\bf{I}_\text{RVE}\in D_{in}}\sum_{i=0}^{n-1}\text{NN}^\text{gate}_i(\bf{I}_\text{RVE})\cdot\log\left(\text{NN}^\text{gate}_{\it i}(\bf{I}_\text{RVE})\right)
\end{equation*}
as introduced by Shannon in \cite{Shannon_Entropy1948}.  Here, D$_\text{in}$ describes the set of inputs for the gating network and $\text{NN}^\text{gate}_i(\bf{I}_\text{RVE})$ is the $i$-th output component for the input vector $\bf{I}_\text{RVE}$. The confidence measurement is then defined as the negative normalized entropy
\begin{equation*}
    confidence = 1.0 - \frac{H}{\log(n)}
\end{equation*}
with the number of experts $n$ as introduced in \Cref{secMoE}. The confidence value increases when most output values of the gating model over the entire dataset are close to zero and it reaches the value $1.0$ when the gating model predicts an output vector with only zeros and a single value of $1.0$. On the other hand, the confidence value decreases for more equally distributed output vectors and reaches zero when the gating model predicts the vector $(\frac{1}{n},\dots,\frac{1}{n})$ for each input.

{\bf Incoherent MoE model with 10 experts.} The first MoE model which we consider is based on ten of the 250 pretrained expert models and a gating network with two hidden layers and 128 neurons in each layer. With regard to the selection of the hyperparameters for the gating architecture, the described model has shown the best results from a wide range of tested models. Given the list of features listed in \Cref{secParam}, the input size of the gating model is equal to ten and the output size is aligned with the number of expert models. In order to define a meaningful selection of expert models that are used for the MoE model, we take the given stress data of the related RVEs into account. All expert models that correspond to one of the 250 RVEs in the training data are considered as potential candidates for the set of experts. We denote the respective expert models by the index related to the respective RVE, which ranges from 0 to 249. The selection of experts is chosen in such a way  such that the stress datasets used to train the respective expert models differ from each other as much as possible. This approach of defining the set of expert models translates to a max-min-dispersion problem \cite{Ravi1994MaxMinDispersion}. We compute the approximate solution of this problem using a greedy algorithm. With this approach of defining the set of experts, we expect to achieve proper generalization with the MoE approach. For the MoE model based on ten experts, the set of expert indices is given by $\lbrace 6, 58, 69, 87, 101, 121, 153, 162, 191, 235\rbrace$.

\begin{figure}
    \begin{minipage}{0.49\textwidth}
        \centering
        \includegraphics[height=4.5cm]{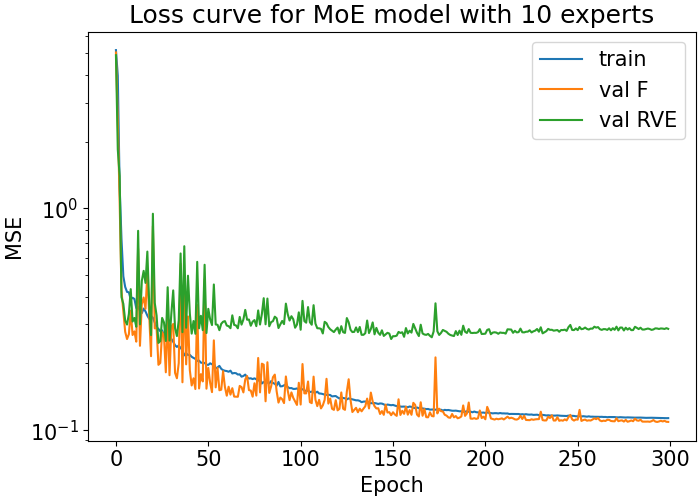}
        \caption{Training loss for the MoE model with ten experts and only gating weights set as trainable. Image and data have first been published in \cite{MagerDiss2026}.}
        \label{MoE60_loss:fig}
    \end{minipage}
    \hspace{0.2cm}
    \begin{minipage}{0.49\textwidth}
        \centering
        \includegraphics[height=4.5cm]{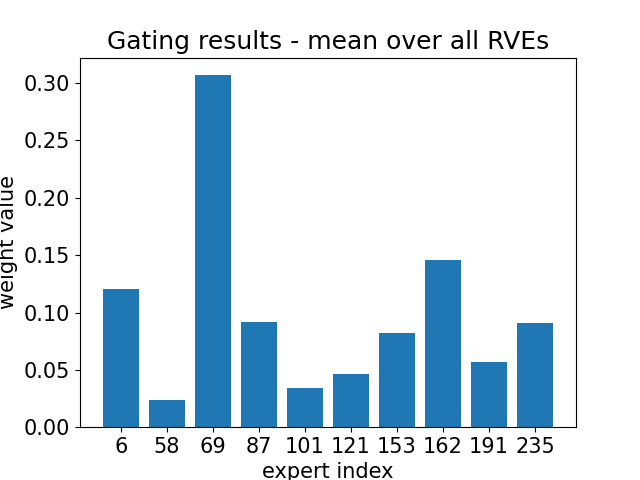}
        \caption{Mean gating output of the MoE model with ten experts and only gating weights set as trainable. Image and data have first been published in \cite{MagerDiss2026}.}
        \label{MoE60_mean_output:fig}
    \end{minipage}
\end{figure}

 The training progress of the incoherently trained MoE model, that is, without changing the expert parameters in phase two, with ten experts is presented in \Cref{MoE60_loss:fig}. In \Cref{MoE60_loss:fig}, the loss curve for the training data as well as for both validation datasets is presented. For the training loss and the $val_{\rm F}$ loss, the most significant decrease in the mean squared error (MSE) is observed. As both loss curves show very similar progression over the scope of 300 training epochs, the given MoE model seems not to be significantly affected by overfitting with regard to the data of the deformation gradients. For the $val_{\rm RVE}$ loss, we see also a significant reduction in the MSE. However, the loss remains relatively constant after about 100 training epochs and does not reduce as significantly as the other two loss values. After training, the $val_{\rm RVE}$ loss equals $2.87e-1$.  The observation that we see only small magnitudes of loss reduction for all three loss terms is largely due to the approach of pre-training the expert models. As the weights of the expert models are already fitted to the stress data of individual RVEs, the MoE model produces already meaningful predictions of the stress components before the gating network is trained. 

 Moreover, we compute the confidence value for the trained gating model based on the input given by a total of 300 RVEs. This includes 250 training RVEs and 50 validation RVEs. The resulting confidence value equals 36\%. Additionally, the weight values computed by the gating model are presented in \Cref{MoE60_mean_output:fig}. The total number of 300 RVEs is used for computing the mean values. In the graph, each bar refers to a single expert model and the height of the bar indicates the mean of the respective weight value. The figure shows that the expert with index 69 has the highest average contribution to the weighted sum that defines the output of the MoE model. Other expert models such as the models with index 58 and 101 have relatively low mean values. Though, all means are distinct from zero which indicates that all expert models provide some contribution to the weighted sum for a given gating input $\bf{I}_{RVE}$.

{\bf Incoherent MoE model with 20 experts.} Since 250 training RVEs are used to create the training data and 50 RVEs are used to create the validation data, a total number of ten expert models that are included in the MoE model appears relatively low. To evaluate whether a higher number of experts would improve the training behavior of the MoE model, we consider a similar approach that includes a total of twenty expert models. To achieve a meaningful comparison between the approaches with different numbers of experts, we use the same incoherent training method. This means that in the second training phase, the weights of the expert models are fixed and the weights of the gating model are adjusted. The indices of the set of experts  
is again computed as the solution of the greedy algorithm applied to the max-min-dispersion problem based on the stress data of the RVEs. For twenty expert models, this results in the set defined by $\lbrace0, 6, 37, 58, 68, 69, 76, 87, 101, 107, 121, 128, 135, 153, 162, 177, 191, 219, 235, 242\rbrace$.

\begin{figure}
    \begin{minipage}{0.49\textwidth}
        \centering
        \includegraphics[height=4.5cm]{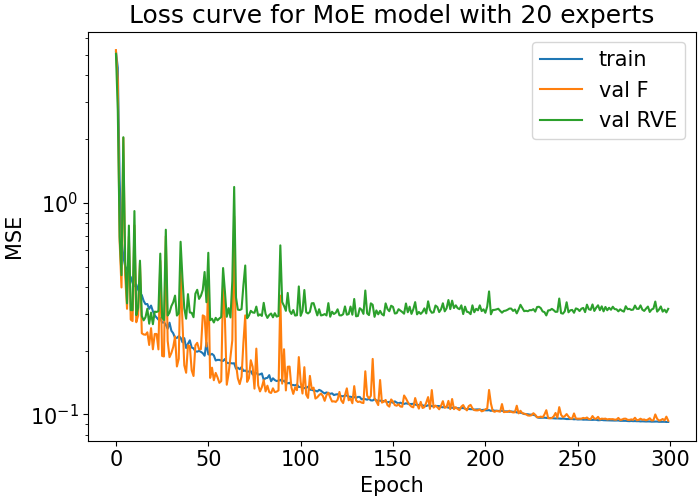}
        \caption{Training loss for the MoE model with twenty experts and only gating weights set as trainable. Image and data have first been published in \cite{MagerDiss2026}.}
        \label{MoE58_loss:fig}
    \end{minipage}
    \hspace{0.2cm}
    \begin{minipage}{0.49\textwidth}
        \centering
        \includegraphics[height=4.5cm]{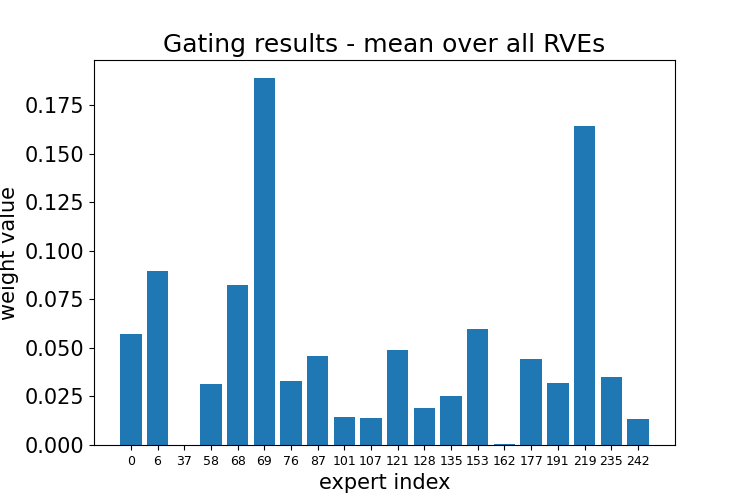}
        \caption{Mean gating output of the MoE model with twenty experts and only gating weights set as trainable. Image and data have first been published in \cite{MagerDiss2026}.}
        \label{MoE58_mean_output:fig}
    \end{minipage}
\end{figure}

For the MoE model with twenty experts, the progression of the loss curves is presented in \Cref{MoE58_loss:fig}. The general training behavior appears very similar to the behavior described for the MoE model with ten experts. The training loss and the validation loss regarding the deformation gradients show close to the same decrease. The $val_{\rm RVE}$ loss curve is not reduced as significantly as the other loss curves. After training of 300 epochs, the $val_{\rm RVE}$ loss equals $3.16e-1$. This MSE value is even higher than the final $val_{\rm RVE}$ loss observed for the MoE model with ten experts. As this is a crucial measurement for the generalizability of the surrogate model with regard to different materials, it is the most important metric for comparing the models. Based on these results, the increase in the number of experts to twenty does not improve the MoE model's performance.

In addition to the evaluation of the training losses, we also analyze the behavior of the gating network for the MoE model with twenty experts. The confidence measurement for the respective gating network equals 40\%. The confidence is slightly higher than the value computed for the MoE model with ten experts. Furthermore, the mean output of the gating model is presented in \Cref{MoE58_mean_output:fig}. The chart shows the highest average weights for the expert models with indices 69 and 219. For the experts with indices 37 and 162, the mean values are approximately zero. This indicates that the respective models have no meaningful contribution to the weighted sum that defines the output of the MoE model. These results indicate that the models in concern could be removed from the set of expert models without an expected change in the model's behavior.

{\bf Coherent MoE models with 10 experts.} In our discussions regarding the previous MoE models, we examined the behavior of the incoherent training approach. In this approach, the weights of the expert models are fixed when the gating model is fitted to the data. In addition to this training method, we investigate the performance of the coherently trained MoE model. This means that the weights of the expert models that result from the first training phase are used to initialize the MoE model. In the second training phase all weights of the MoE model including the expert weights and weights of the gating model are set trainable. To minimize the MSE loss, all weights of the MoE model are adjusted to fit the training data. As the MoE model with ten experts shows better training performance for the incoherent training approach, we consider the same model architecture for the coherent training approach. Furthermore, the same selection of expert models is used as a basis for the MoE model.

\begin{figure}
    \begin{minipage}{0.49\textwidth}
        \centering
        \includegraphics[height=4.5cm]{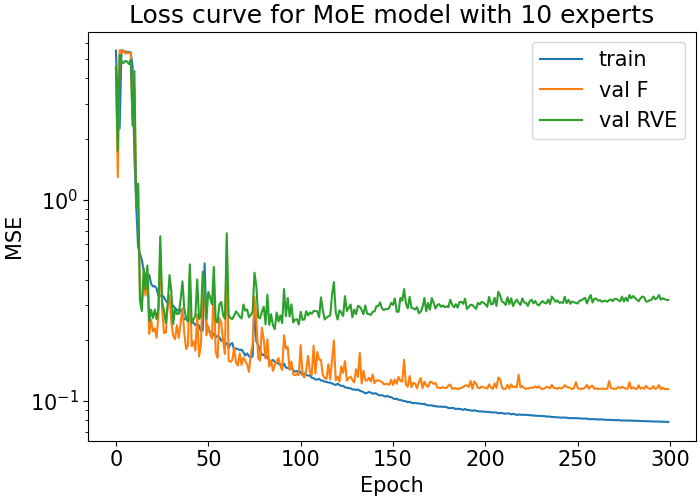}
        \caption{Training loss for the ten experts MoE model with gating and expert weights set as trainable. Image and data have first been published in \cite{MagerDiss2026}.}
        \label{MoE61_loss:fig}
    \end{minipage}
    \hspace{0.2cm}
    \begin{minipage}{0.49\textwidth}
        \centering
        \includegraphics[height=4.5cm]{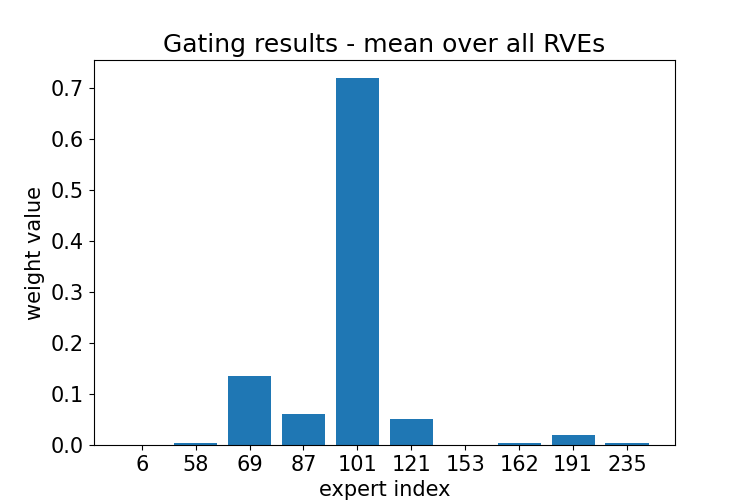}
        \caption{Mean gating output of the ten experts MoE model with gating and expert weights set as trainable. Image and data have first been published in \cite{MagerDiss2026}.}
        \label{MoE61_mean_output:fig}
    \end{minipage}
\end{figure}

The progression of the loss values for the coherent training method is presented in \Cref{MoE61_loss:fig}. With this approach, we see a slightly better reduction of the training loss compared to the incoherently trained models. However, the difference between the training loss and both validation losses is significantly larger than in \Cref{MoE60_loss:fig,MoE58_loss:fig}. The MSE regarding the $val_{\rm RVE}$ data equals $3.17e-1$ after 300 epochs of training. The gaps between training and validation loss curves indicate that the coherent training approach leads to a more notable overfitting with respect to the training data. Based on these results, adjusting the expert weights in the second training phase does not improve the training performance of the MoE model.

Additionally, the mean output of the gating model resulting from the coherent training approach is presented in \Cref{MoE61_mean_output:fig}. The average weight values are computed based on all 300 RVEs. The chart shows that the expert with index 101 has by far the highest average weight. With a mean weight of approximately 0.7, the output of this expert model has the highest average contribution to the weighted sum that defines the output of the MoE model. The other expert models have only minor contributions to the weighted sum. This focus of the gating model on the expert with index 101 leads also to a higher confidence value of 60\%.

We see a very similar behavior for other coherently trained MoE models which consider a higher or lower number of experts. Although the models yield lower training loss values, the $val_{\rm RVE}$ loss ends up higher than in the respective incoherent case. The training time of the coherent models is also significantly higher due to the larger number of weights which are optimized.

\subsection{Importance of the input features}
In this section, we investigate the importance of the ten different input features of the gating network related to pore sizes, beam orientations, and the volume fraction of the RVE. We consider the permutation feature importance \cite{Altmann2010PermutationImportance} as a measurement. The concept of this measurement is to randomly permute a single input feature of a given dataset and compare the machine learning model output obtained for the modified dataset with the corresponding output for the original dataset. To define the permutation feature importance, we consider the input dataset D$_\text{in}$ and, for a given input sample $\bf{I}_\text{RVE}$, introduce the permuted sample $\bf{\hat{I}}_{\text{RVE},j\leftarrow \pi}$. This sample coincides with $\bf{I}_\text{RVE}$ for all entries except feature $j$, whose value is randomly exchanged with the $j$-th feature values of the other samples in the dataset. We define $\Psi\in\R^{n\times n_{in}}$ as the matrix collecting the importance measures for each feature and each output. Accordingly, the entry $\Psi_{i,j}$, representing the importance of feature $j$ with respect to output $i$, is defined as
\begin{equation*}
    \Psi_{i,j} = \frac{1}{\lvert \text{D}_\text{in}\rvert}\sum_{\bf{I}_\text{RVE}\in D_{in}}\left(\text{NN}^\text{gate}_i({\bf I}_\text{RVE}) - \text{NN}^\text{gate}_i({\hat{\bf I}}_{\text{RVE},j\leftarrow \pi})\right)^2
\end{equation*}
for a given gating network $\text{NN}^\text{gate}$.

\begin{figure}
    \centering
    \includegraphics[height=6.5cm]{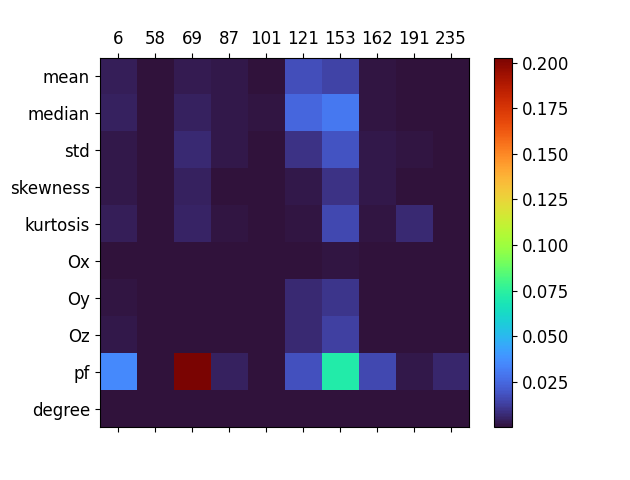}
    \caption{Permutation feature importance of the gating network with ten experts trained with the weights of the expert models fixed. Image and data have first been published in \cite{MagerDiss2026}.}
    \label{MoE60_importance:fig}
\end{figure}

The importance measurements for the gating model which uses ten experts and which is trained with the fixed weights of the expert models is shown in \Cref{MoE60_importance:fig}. Each field of the colormap represents one combination of a feature and an output value. The color of the respective field illustrates the respective importance value. Here, higher values indicate a stronger importance of the feature for the output. The figure implies that the output weights of the gating model have the highest sensitivity to the porous fraction.  The highest importance value is computed for the relationship between the porous fraction and the output of expert 69. The high importance of the porous fraction for the mechanical properties of an open-porous material is also investigated in experimental analysis for biopolymer aerogels \cite{DLR_biopoly_aerogel2016,DLR_biopoly_aerogel2018}. In addition to the porous fraction, the measurements referring to the pore-size distribution show also meaningful importance for the output of the gating model.

The measurements for different gating models show very similar results and also indicate the highest importance with the porous fraction.

\section{Application of the mixture of experts surrogate model in multiscale simulations}
The main purpose of our MoE approach is the online application of the model in multiscale simulations. The MoE surrogate model features the advantage that the simulation of different materials is possible with a single model that has exclusively been trained for an available set of materials.  As mentioned before, in the previous NN approach in \cite{AerogelFE2}, a new surrogate model is required to be trained for every material considered for a simulation. In this section, we discuss the performance of the novel MoE approach when applied as a surrogate model on the microscopic scale to the FE$^2$ method. For this purpose, we compare the macroscopic results using three different microscopic solvers: the MoE approach, the basic NN approach as described in \cite{AerogelFE2}, and the consecutive beam frame solve.  The latter one can be seen as the ground truth.

\subsection{Macroscopic results in comparison to beam frame approach}
To evaluate the macroscopic performance of the MoE approach for practical applications, we consider the torsion of a cube geometry and we compare the FE$^2$ solutions computed with the MoE model or the beam frame solver on the microscale. The macroscopic mesh is discretized using $20\,480$   $\mathcal{P}_1$, that is, linear, finite elements with $14\,739$ degrees of freedom (dof). We compute a torsion of 18 degree relative to the $x$-axis. Dirichlet boundary conditions are applied to the two faces of the cube pointing in $x$-direction.

For the microscopic problem, we consider an RVE that has not been included in the training data for the MoE model, that is, no explicit expert model has been included for this specific RVE. The beam frame structure contains $4\,870$ beam elements connected at $2\,988$ vertices. For the beam frame solver, this results in a system of equations with $17\,928$ dof.

\begin{figure}
    \centering
    \includegraphics[height=5cm]{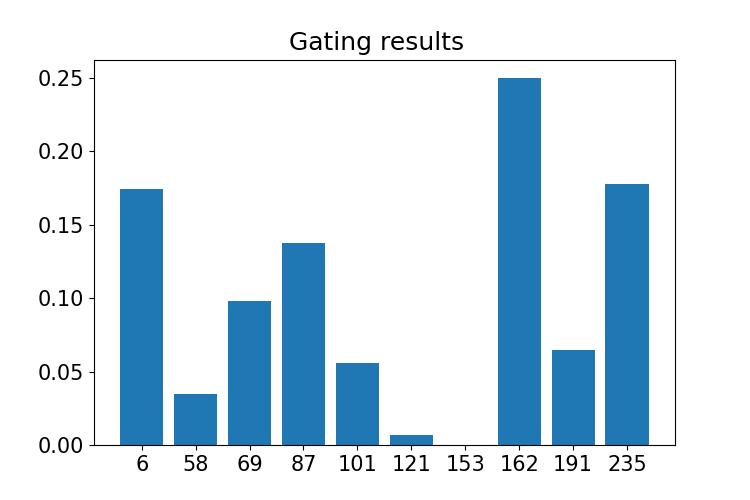}
    \caption{Computed output of the gating network with ten experts for the considered input RVE.}
    \label{MoE60_gating274:fig}
\end{figure}

For solving the FE$^2$ problem with the MoE approach, we consider the MoE model with ten experts which was trained incoherently in two stages, since it provides the best performance regarding the $val_{\rm RVE}$ loss within the studied framework. 
Using the parameterized representation of the example beam frame RVE, the gating network yields an output vector with ten nonzero values. The output values are shown in \Cref{MoE60_gating274:fig} with the respective expert indices on the $x$-axis. For the considered feature vector $\textbf{I}_{RVE}$, the experts with indices 6, 162, and 235 contribute the most to the weighted sum that defines the output of the MoE model. Note that the output of the expert models with indices 121 and 153 are not computed in the evaluation of the MoE model within the FE$^2$ method, as the weights for these models are close to zero. Neglecting the evaluation of models with corresponding weight values close to zero reduces the computational effort for the evaluation of the surrogate model. We use the threshold of $w_{min}=\frac{0.1}{n}$ that depends on the respective number of experts given within the MoE architecture. The weights of the remaining experts are rescaled to ensure the sum over all weights to be equal to one.

We use Newton's method as a nonlinear solver at the macroscopic level with both the beam frame model and the MoE approach. However, in case of the beam frame model, the macroscopic Jacobian is not computed exactly but an approximation with central difference quotients is used; see~\cite{AerogelFE2} for details. For the MoE model, the Jacobian can be computed directly as the weighted sum of the Jacobians of the single expert NNs.

\begin{figure}
    \centering
    \includegraphics[height=7cm, trim={0cm 0cm 0cm 0cm}, clip]{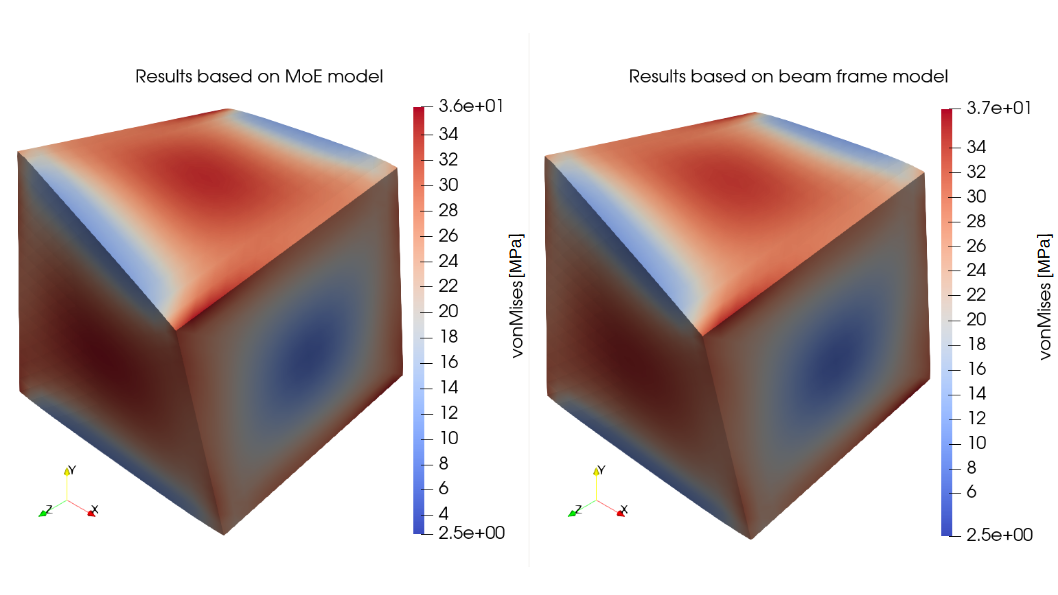}
    \caption{Results for the torsion of a cube geometry computed with the MoE model (left) and with the beam frame model (right) on the microscale of an FE$^2$ simulation. The colors in the plots represent the von Mises stress. Images and data have first been published in \cite{MagerDiss2026}.}
    \label{macro_results251:fig}
\end{figure}

\begin{figure}
    \centering
    \includegraphics[height=7cm, trim={0cm 0cm 0cm 0cm}, clip]{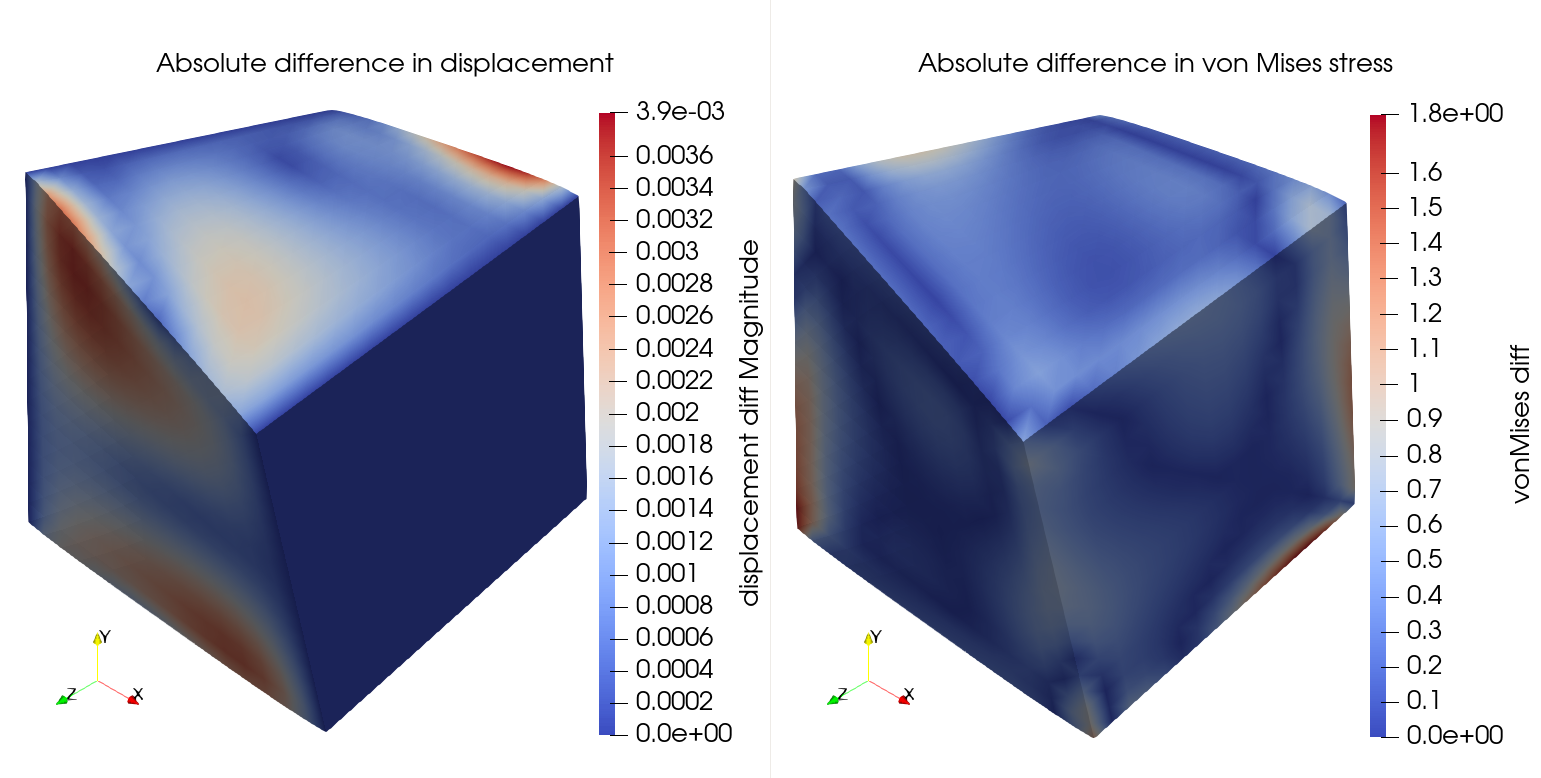}
    \caption{Deviation between the FE$^2$ solution using the MoE surrogate model and the beam frame solver 
    on the microscale. Left image shows the deviation regarding the deformation with the colors representing the norm over the deviation vector in each examination point. Right image shows the deviation regarding the von Mises stress with the colors representing the absolute difference of the stress measurements in each point. Images and data have first been published in \cite{MagerDiss2026}.}
    \label{macro_diff251:fig}
\end{figure}

The results of both simulations are shown in \Cref{macro_results251:fig} with the colors representing the von Mises stress values. Visually, there are no significant differences observable between the two solutions. To emphasize the differences between both solutions, the plot on the left in \Cref{macro_diff251:fig} highlights the difference in the deformation between the beam frame solution and the approximated simulation computed with the MoE model. A similar illustration of the difference in the computed von Mises stress is presented on the image on the right. The colors represent the absolute difference between the stress measurements in each Gauss evaluation point. The error plots show some regions within the geometry with a higher deviation between the solutions, but, overall, the maximum deviation of $3.9\cdot 10^{-3}$ between the deformation values is relatively small. Note that the difference at the opposing surfaces in $x$-direction is zero, as the same Dirichlet boundary conditions are applied to the macroscopic problems. For the measurements of the von Mises stress, the absolute error is relatively evenly distributed. 

Over the total geometry, the relative deviations regarding the deformation and the von Mises stress are given by
\begin{align*}
    \frac{\lVert\overline{u}_{BF}-\overline{u}_{MoE}\rVert}{\lVert\overline{u}_{BF}\rVert}=&1.2244\cdot 10^{-2},  \\
    \frac{\lVert\overline{\sigma}_{vM,BF}-\overline{\sigma}_{vM,MoE}\rVert}{\lVert\overline{\sigma}_{vM,BF}\rVert}=&1.1937\cdot 10^{-2}
\end{align*}
with the macroscopic deformation vector resulting from the beam frame solution denoted by $\overline{u}_{BF}$ and the deformation resulting from the MoE solution denoted by $\overline{u}_{MoE}$. Similarly, the von Mises stresses are given by $\overline{\sigma}_{vM,BF}$ and $\overline{\sigma}_{vM,MoE}$.

\subsection{Comparison with the surrogate approach based on a single neural network}
The expert models within the MoE approach are based on the machine learning-based surrogate model introduced in \cite{AerogelFE2}. Since the MoE model can be considered as a further development of the consecutive NN approach, it is relevant to compare the performance of both models for different RVEs. The MoE model yields the advantage that the gating model  explicitly uses statistical parameters of a specific RVE as input data to compute different outputs for different RVEs with the same deformation gradient. The basic NN approach from~\cite{AerogelFE2}, on the other hand, uses exclusively the deformation gradient as input. However, since the RVEs considered in our dataset are all based on the nanostructure of cellulose aerogels with pore-size distributions in the same range, we expect the macroscopic deformations for different RVEs from our dataset to be similar. Therefore, in this section, we investigate the differences between the macroscopic results computed with the MoE approach and the results obtained using only the NN of a single expert model. The results presented in \Cref{MoE60_mean_output:fig,MoE58_mean_output:fig} show that the expert neural network with index 69 has the highest average contribution to the weighted sum of the incoherently trained MoE models. Based on this performance, it is expected that this expert NN provides a meaningful basis for the prediction of various RVEs within the training and validation datasets. To have a proper expert NN for the comparison, we use the expert NN with index 69 for the comparison with the MoE approach. Note again that this expert NN has only been trained for a single RVE and since the only input is the deformation gradient, it cannot adapt the mechanical behavior for different RVEs.

\begin{figure}
    \begin{minipage}{0.5\textwidth}
        \includegraphics[width=\textwidth]{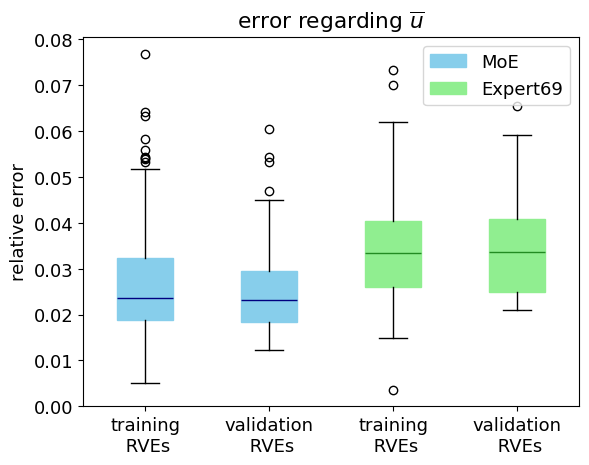}
    \end{minipage}
    \begin{minipage}{0.5\textwidth}
        \includegraphics[width=\textwidth]{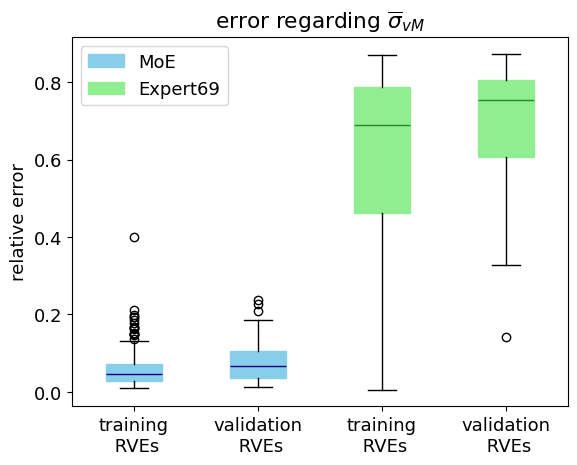}
    \end{minipage}
    \caption{Distribution of the error between the surrogate model and the beam frame model for different RVEs. In the figures, the error is evaluated regarding the deformation (left) and the von Mises stress (right). Blue boxplots show the error regarding the MoE model and green boxplots show the error regarding the expert NN with index 69. Images and data have first been published in \cite{MagerDiss2026}.}
     \label{boxplots:fig}
\end{figure}

The comprehensive results of the comparison are presented in \Cref{boxplots:fig}. The boxplot diagram on the left shows the distribution of the error with regard to the deformation variable and similarly, the figure on the right illustrates the deviation with regard to the resulting von Mises stress values. In each of the diagrams, the blue boxes represent the differences between the MoE model and the beam frame models. We further differentiate the error distributions for the training and the validation RVEs. For the boxplot referring to the training RVEs, we compute the $L_2$-error between the two computed macroscopic deformation vectors. This results in 250 error values, the distribution of which is shown within the boxplot. Similarly, the box referring to the validation RVEs presents the distribution of 50 error values regarding the same computations for the validation RVEs. For the results computed with the single NN of expert 69, the distribution of the errors are illustrated with the green boxes in \Cref{boxplots:fig}. For all presented boxplots, the length of the whiskers is computed based on the 1.5 interquartile range (IQR) value.

Regarding the boxplot that visualizes the error in deformation for expert 69, note that the single outlier in the error dataset for the training RVEs with a very low error value close to 0 belongs to the RVE with index 69 the NN is explicitly trained for. With regard to the deformation $\overline{u}$, the results in \Cref{boxplots:fig} indicate that the MoE approach largely yields lower error values than the approach of using only the single NN on the microscopic scale. While the difference between the distributions does not seem to be great, the lower errors for the MoE model are apparent for the training RVEs as well as for the validation RVEs. The mean error values for the MoE model and the NN trained for a single RVE are shown in \Cref{mean_error:tab}. The mean values underscore the model's improved performance, which was already apparent from \Cref{boxplots:fig}.

 The difference in the error regarding the macroscopic deformation variable $\overline{u}$ is relatively small between the results based on the MoE model and on the single NN, as differences in the microstructure cause only small deviations in the macroscopic deformation. The influence of the microstructure on the macroscopic stress value is more meaningful, since the amount of force required for deformation varies depending on the pore composition and the thickness of the fibrils. Therefore, a more significant difference between the MoE approach and the results based on expert 69 can be observed with regard to the error computed for the von Mises stress. In the diagram on the right in \Cref{boxplots:fig}, the errors of the MoE model are considerably smaller in comparison to the errors resulting from the single NN of expert 69. This finding is supported by the mean values in \Cref{mean_error:tab} which show that the error of the MoE model is on average about one magnitude smaller than the error of the single NN. This significant difference is largely driven by the difference in scale of the stress values resulting from different RVEs. While the stress values predicted by expert 69 are always in the same range, the MoE model is capable of adapting to the individual stress scales. 

Regarding the deviation in the deformation value $\overline{u}$ for the training RVEs, the MoE and the single NN approach both show a significant number of outliers with errors above the end of the upper whisker. For the MoE model, the dataset features eight outliers and for the NN model, two outliers are present within the dataset.  In sum, the relative number of outliers is rather small regarding the total number of 250 RVEs for this dataset. However, especially for the MoE approach, it is desired that the method provides small errors for every considered RVE within the dataset. Reducing the error for all RVEs and also reducing the number of outliers poses an important challenge for the future work regarding further optimization of the MoE approach.

\begin{table}
	\begin{tabular}{|l|cc|cc|}
        \hline
        & \multicolumn{2}{|c|}{$\overline{u}$} & \multicolumn{2}{|c|}{$\overline{\sigma}_{vM}$} \\

		& MoE model & NN for expert 69 & MoE model & NN for expert 69 	\\
		\hline
		Training RVEs &$2.63\cdot 10^{-2}$  & $3.42\cdot 10^{-2}$ &$5.78\cdot 10^{-2}$  & $6.20\cdot 10^{-1}$  \\
		Validation RVEs &$2.58\cdot 10^{-2}$  & $3.40\cdot 10^{-2}$ &$8.05\cdot 10^{-2}$  & $4.35\cdot 10^{-1}$  \\
		\hline
		Total &$2.62\cdot 10^{-2}$  & $3.42\cdot 10^{-2}$ &$6.16\cdot 10^{-2}$  & $6.32\cdot 10^{-1}$  \\
		\hline
	\end{tabular}
	\caption{Mean error regarding the macroscopic deformation and von Mises stress for the MoE approach and the single expert NN with index 69. Data has first been published in \cite{MagerDiss2026}.}
	\label{mean_error:tab}
\end{table}

In addition to the described computations with the single NN of expert 69, we also investigated the performance of the MoE compared to the computations with other expert NN, each trained for a single RVE. The basic NN approaches showed overall very similar behavior as the NN for expert 69 with the MoE model outperforming each of the considered models. Since adding the results for further models would not have provided any further insight, we only included the results for expert 69.

\section{Conclusion}\label{secConclusion}
The MoE model offers a useful extension to the multiscale approach introduced in \cite{AerogelFE2}. The model shows good results in terms of generalizability, as the error regarding materials not included in the training remains relatively small. It reduces the computational effort for comparing simulation results regarding different open-porous materials since it eliminates the need of additional training for each given RVE.

The comparison of different MoE architectures has shown that the number of experts plays a significant role for the models performance. A higher number of experts can lead to overfitting while a number that is too low leads to worse training losses. For our applications, the MoE model with ten experts provides the best trade-off between predictive performance and model complexity.   

 For evaluating the macroscopic performance of the MoE approach, we compared the model to the beam frame solver as well as to the basic NN approach which uses a single feedforward NN. The results exhibit that the MoE solver yields a performance similar to that of the beam frame solver and also reduces the error in comparison to the basic NN approach. 
 However, if only the deformation behavior is considered, the difference between the error of the MoE model and the single feedforward NN is relatively small due to the similarity of the considered RVEs. With regard to the error in the computation of stress measurements, the MoE model shows greater advantages over the single NN approach.

 Applying the MoE approach to a broader range of porous materials could further highlight the model’s strengths compared to simple NN approaches. Previous work has shown that training a gating network within the MoE model benefits when the individual experts are negatively correlated \cite{Masoudnia2014MoE_Survey}. In our application to fibrillar nanostructured materials and aerogels, the individual experts are strongly correlated. With materials that differ more significantly, one might expect a weaker correlation and, consequently, improved training performance of the gating model. We view the results presented in this work as a demonstration of the MoE approach, which can be applied to other applications.

For future work, the architecture of the gating model could be further examined by using the direct representation of the beam frame model as graph, the embedding of graph neural networks (GNN) as gating models would be possible. Further investigation is required to evaluate the potential of this approach.

\backmatter
\bmhead{Acknowledgements}
The authors gratefully acknowledge the use of the computational
facilities of the Center for Data and Simulation Science (CDS) of the University of
Cologne.
\bmhead{Author contributions}
All authors are listed in alphabetical order. The contributions are specified below using the CRediT taxonomy.
AK: Conceptualization, Methodology, Supervision, Writing - review \& editing. ML: Methodology, Visualization, Writing - review \& editing. 
LM: Data Curation, Investigation, Methodology, Software, Validation, Visualization, Writing - original draft. AR: Methodology, Writing - review \& editing. JWH: Methodology, Writing - review \& editing


\bibliography{sn-bibliography}

\end{document}